\documentclass[letterpaper,10pt,conference]{ieeeconf}
\IEEEoverridecommandlockouts 
\usepackage{amssymb,amsmath,graphicx,subfigure,warmread,times}
\usepackage[all,import]{xy}

\newtheorem{theorem}{Theorem}

\newcommand{\norm}[1]{\ensuremath{\left\| #1 \right\|}}
\newcommand{\bracket}[1]{\ensuremath{\left[ #1 \right]}}
\newcommand{\braces}[1]{\ensuremath{\left\{ #1 \right\}}}
\newcommand{\parenth}[1]{\ensuremath{\left( #1 \right)}}
\newcommand{\refeqn}[1]{(\ref{eqn:#1})}
\newcommand{\reffig}[1]{Fig. \ref{fig:#1}}
\newcommand{\tr}[1]{\mbox{tr}\ensuremath{\bracket{#1}}}
\newcommand{\deriv}[2]{\ensuremath{\frac{\partial #1}{\partial #2}}}
\newcommand{\SO}{\ensuremath{\mathrm{SO(3)}}}
\newcommand{\T}{\ensuremath{\mathrm{T}}}
\newcommand{\so}{\ensuremath{\mathfrak{so}(3)}}

\renewcommand{\Re}{\ensuremath{\mathbb{R}}}
\renewcommand{\S}{\ensuremath{\mathbb{S}}}

\newcommand{\wh}{\widehat}
\newcommand{\wt}{\widetilde}
\newcommand{\la}{\label}
\newcommand{\be}{\begin{equation}}
\newcommand{\ee}{\end{equation}}
\newcommand{\bea}{\begin{eqnarray}}
\newcommand{\eea}{\end{eqnarray}}
\newcommand{\beas}{\begin{eqnarray*}}
\newcommand{\eeas}{\end{eqnarray*}}
\newcommand{\nn}{\nonumber}
\newcommand{\cJ}{\mathcal{J}}
\newcommand{\Tp}{^{\mbox{\small T}}}
\date{}

\renewcommand{\thesubsection}{\arabic{subsection}. }
\renewcommand{\thesubsubsection}{\arabic{subsection}.\arabic{subsubsection} }

\begin{document}
\title{Global Attitude Estimation using Uncertainty Ellipsoids}
\author{Taeyoung Lee$^1$,\ Amit K. Sanyal$^2$,\ Melvin Leok$^3$,\ N. Harris 
McClamroch$^1$ \\
\parbox{3 in}{\centering $^1$Department of Aerospace Engineering, \\
$^3$Department of Mathematics, \\
       University of Michigan, \\
        Ann Arbor, MI 48109  \\
	{\tt \{tylee, mleok, nhm\}@umich.edu }} 
\hspace*{-0.5in}
\parbox{4 in}{\centering $^2$Department of Mechanical and \\ Aerospace Engineering, \\
Arizona State University, \\
Tempe, AZ 85287 \\
{\tt sanyal@asu.edu}}
}

\maketitle
\begin{abstract}
Attitude estimation is often a prerequisite for control of the attitude or 
orientation of mechanical systems. Current attitude estimation algorithms use 
coordinate representations for the group of rigid body orientations. All 
coordinate representations of the group of orientations have associated problems. 
While minimal coordinate representations exhibit kinematic singularities for 
large rotations, non-minimal coordinates like quaternions require satisfaction 
of extra constraints. A deterministic attitude estimation problem for a rigid 
body with bounded measurement errors is considered here. An attitude estimation 
algorithm that globally minimizes the attitide estimation error, is obtained. 
Assuming that the initial attitude, the initial angular velocity and measurement 
noise lie within given ellipsoidal bounds, an uncertainty ellipsoid that
bounds the attitude and the angular velocity of the rigid body is
obtained. The center of the uncertainty ellipsoid provides point
estimates, and the size of the uncertainty ellipsoid measures the
accuracy of the estimates. The point estimates, and the uncertainty ellipsoids are propagated using a Lie group variational integrator, and its linearization, respectively. The attitude estimation is optimal in the
sense that the attitude estimation error and the size of the uncertainty 
ellipsoid is minimized.
\end{abstract}

\section{Introduction}

Attitude estimation is often a prerequisite for controlling aerospace and 
underwater vehicles, mobile robots, and other mechanical systems moving in 
space. Hence, attitude estimation may be used in spacecraft and aircraft, 
unmanned vehicles and robots, including walking robots. In this paper, we look 
at the attitude estimation problem for the uncontrolled dynamics of a rigid 
body in an attitude-dependent force potential. The estimation scheme we present 
has the following important features: (1) the attitude is globally represented 
without using any coordinate system, (2) the filter obtained is not a Kalman or 
extended Kalman filter, and (3) the attitude and angular velocity measurement 
errors are assumed to be bounded, with ellipsoidal uncertainty bounds. The 
static attitude estimation using a global attitude representation is based 
on~\cite{san2006}. Such a global representation has been recently used for 
partial attitude estimation with a linear dynamics model in~\cite{rehu}. 

The attitude determination problem for a rigid body from vector measurements was 
first posed in \cite{wah}. A sample of the literature in spacecraft attitude 
estimation can be found in \cite{bosh, cram, mark, shus1, shus2}. Applications of 
attitude estimation to unmanned vehicles and robots can be found in 
\cite{rehu, badu,rosb, vafo}. Most existing attitude estimation schemes use 
coordinate representations of the attitude. As is well known, minimal coordinate 
representations of the rotation group, like Euler angles, Rodrigues parameters, 
and modified Rodrigues parameters (see \cite{cram2}), usually lead to geometric or 
kinematic singularities. Non-minimal coordinate representations, like the 
quaternions used in the quaternion estimation (QUEST) algorithm and its several 
variants (\cite{bosh, shus2, psia}), have their own associated problems. Besides 
the extra constraint of unit norm that one needs to impose on the quaternion, the 
quaternion representation for a given rotation depends on the sense of rotation 
used to define the principal angle, and hence can be defined in one of two ways. 

A brief outline of this paper is given here. In Section II, the attitude 
determination problem for vector measurements with measurement noise is 
introduced, and a global attitude determination algorithm which minimizes 
the attitude estimation error is presented. In Section III, the attitude dynamics 
and dynamic estimation problem is formulated, and an algorithm to numerically 
integrate the dynamics is presented. Section IV presents the attitude 
estimation scheme with attitude and angular velocity measurements. Section V 
presents some simulation results followed by conclusions in Section VI.

\section{Attitude Determination}

\subsection{Attitude determination from vector observations}

Attitude of a rigid body is defined by the orientation of a body
fixed frame with respect to a reference frame, and the attitude is
represented by a rotation matrix that is a $3\times 3$ orthogonal
matrix with determinant of 1. Rotation matrices have a group
structure denoted by $\SO$. The group operation of $\SO$ is matrix
multiplication, and its action on $\Re^3$ takes a vector represented 
in body fixed frame into the reference frame by matrix multiplication.

We denote the known direction vector of the $i$th point in the reference 
frame as $e^i\in\S^2$, and the corresponding vector represented in the
body fixed frame as $b^i\in\S^2$. These direction vectors are
normalized so that they have unit lengths. The $e^i$ and $b^i$ are 
related by a rotation matrix $C\in\SO$ that defines the rigid body 
attitude;
\begin{align*}
e^i = C b^i,
\end{align*}
for all $i\in\braces{1,2,\cdots,m}$, where $m$ is the number of
measurements. We assume that $e^i$ is known accurately and $b^i$ is 
measured by sensors in the body fixed frame. Let the measured direction 
vectors (with sensor errors) be denoted $\tilde{b}^i\in\S^2$, and let 
an estimate of the rotation matrix be denoted $\wh{C}\in\SO$. The 
estimation error is given by
\begin{align*}
e^i-\wh{C}\tilde{b}^i.
\end{align*}

The attitude determination problem consists of finding an estimate 
$\wh{C}\in\SO$, and is given by the following weighted least squares 
problem:
\begin{align}
\min_{\wh{C}} \mathcal J & =\frac{1}{2}\sum_{i=1}^m w_i
(e^i-\wh{C}\tilde b^i)\Tp (e^i-\wh{C}\tilde b^i),\nonumber\\
& = \frac{1}{2} \tr {(E-\wh{C}\tilde B)\Tp W(E-\wh{C}\tilde B)},\label{eqn:wahba}\\
&\text{subject to } \wh{C}\in\SO,\nn
\end{align}
where $E=\bracket{e^1,e^2,\cdots,e^m}\in\Re^{3\times m}$, $\tilde
B=\bracket{\tilde b^1,\tilde b^2,\cdots,\tilde b^m}\in\Re^{3\times
m}$, and $W=\mathrm{diag}\bracket{w^1,w^2,\cdots,w^m}\in\Re^{m\times
m}$ has weight factors for each measured vector. 

This problem is known as Wahba's problem~\cite{wah}. The solution in 
terms of quaternions, known as the QUEST algorithm, is presented 
in~\cite{shus1}. A solution without using generalized attitude coordinates 
is given in~\cite{san2006}. A necessary condition for optimality of 
\refeqn{wahba} is given by
\begin{align}
L\Tp \hat{C}=\hat{C}\Tp L,\label{eqn:san2006}
\end{align}
where $L=EW\wt B\Tp \in\Re^{3\times 3}$.  

The following result, which is proved in~\cite{san2006}, gives an 
unique estimate $\wh{C}\in\SO$ of the attitude matrix that 
solves the attitude determination problem \refeqn{wahba}.  
\begin{theorem}
The unique minimizing solution to the attitude determination problem 
\refeqn{wahba} is given by 
\be \wh{C}= SL,\;\ S=Q\sqrt{(RR\Tp)^{-1}}Q\Tp, \la{soln} \ee
where 
\be L=QR,\;\ Q\in\SO, \la{LQR} \ee
and $R$ is upper triangular and invertible; this is the QR decomposition 
of $L$. The symmetric positive definite (principal) square root is used 
in (\ref{soln}).
\la{atde}
\end{theorem}
The proof is based on the fact that $\cJ$ is a Morse function, i.e., its critical 
points are non-degenerate. 
From the Morse lemma \cite{miln}, we conclude that these non-degenerate critical 
points are isolated, and hence the estimate given by (\ref{soln}) uniquely 
minimizes the attitude estimation error. 

\subsection{Estimation with bounded state uncertainties}

A stochastic state estimator requires probabilistic models for 
the state uncertainty and the noise, which are often not available. 
Assumptions are usually made on the statistics of disturbance and noise 
processes, in order to make the estimation problem mathematically 
tractable. In many practical situations such idealized assumptions are 
not appropriate, and may cause poor estimation performance~\cite{TheSkaSou}.
An alternative deterministic approach is to specify bounds on the
uncertainty and the measurement noise without any assumptions on their
distribution. Noise bounds are available in many cases, and such a 
deterministic estimation scheme is robust to the noise distribution. 
An efficient but flexible way to describe the bounds is using
ellipsoidal sets, referred to as uncertainty ellipsoids.

\renewcommand{\xyWARMinclude}[1]{\includegraphics[width=0.45\textwidth]{#1}}
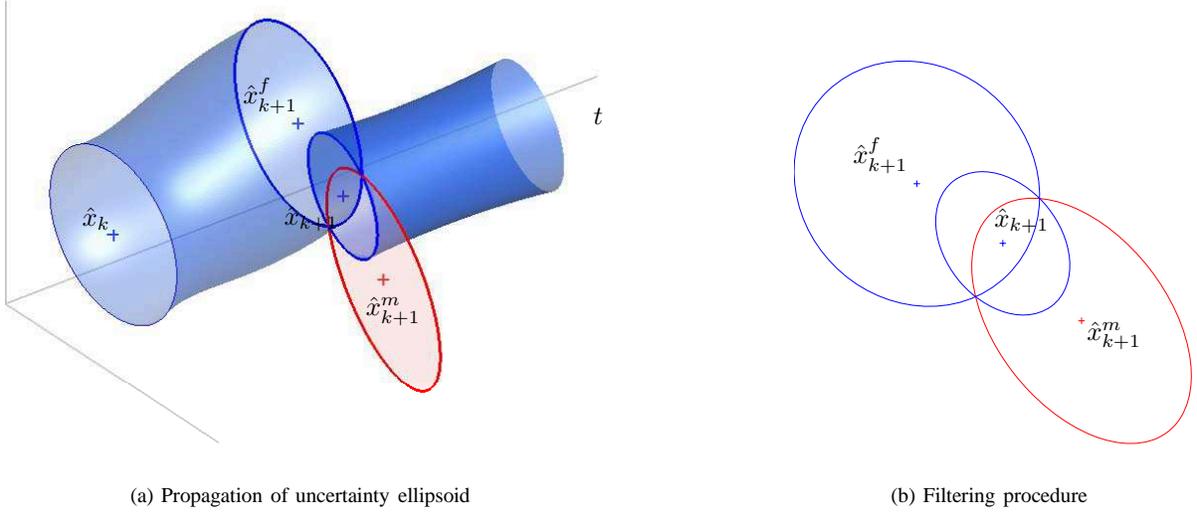
\begin{figure*}[t]
    \centerline{\subfigure[Propagation of uncertainty ellipsoid]{
    $$\begin{xy}
    \xyWARMprocessEPS{tube}{eps}
    \xyMarkedImport{}
    \xyMarkedMathPoints{1-5}
    \end{xy}$$}
    \hspace*{2cm}
    \renewcommand{\xyWARMinclude}[1]{\includegraphics[width=0.3\textwidth]{#1}}
    \subfigure[Filtering procedure]{
        $$\begin{xy}
    \xyWARMprocessEPS{tubesection}{eps}
    \xyMarkedImport{}
    \xyMarkedMathPoints{1-3}
    \end{xy}$$}}
    \caption{Uncertainty ellipsoids}\label{fig:ue}
\end{figure*}

This deterministic estimation procedure for a 2 dimensional system
is illustrated in \reffig{ue}, where the left figure shows time 
evolution of an uncertainty ellipsoid, and the right figure shows a
cross section at a fixed time when the state is measured. At the
$k$th time step, the state is bounded by an uncertainty ellipsoid
centered at $\hat{x}_k$. This initial ellipsoid evolves over time. 
Depending on the dynamics of the system, the size and the shape of 
the tube are changed. At the $k+1$th time step, the predicted uncertainty 
ellipsoid is centered at $\hat{x}_{k+1}^f$. The state is then measured by 
sensors, and another ellipsoidal bound on the state is obtained by the 
measurements. The measured uncertainty ellipsoid is centered at 
$\hat{x}_{k+1}^m$. The state lies in the intersection of the two ellipsoids. 
In the estimation procedure, we find a new ellipsoid that contains the
intersection, which is shown in the right figure. The center of the
new ellipsoid, $\hat{x}_{k+1}$ is considered as a point estimate at
time step $k+1$, and the magnitude of the new uncertainty
ellipsoid measures the accuracy of the estimation. 
This deterministic estimation is optimal in the sense that the size of 
the new ellipsoid is minimized.

A deterministic estimation process based on set theoretic results was 
developed in~\cite{Sc68}. Optimal deterministic estimation is 
considered in~\cite{MaNo96} and \cite{DuWaPo01}, where an analytic 
solution for the minimum ellipsoid that contains a union or an 
intersection of ellipsoids is obtained.

\section{Attitude Dynamics and Dynamic Attitude Estimation}

\subsection{Equations of motion}
We now consider dynamic state estimation of the attitude dynamics of 
a rigid body in a potential $U(C):\SO\mapsto\Re$ determined by 
the attitude, $C\in\SO$. A spacecraft on a circular orbit including 
gravity gradient effects~\cite{acc06}, or a 3D pendulum~\cite{cca05} 
can be so modeled. The continuous equations of motion are given by
\begin{gather}
J\dot\omega + \omega\times J\omega = M,\\
\dot{C} = C S(\omega),\label{eqn:Rdot}
\end{gather}
where $J\in\Re^{3\times 3}$ is the moment of inertia matrix of the
rigid body, $\omega\in\Re^3$ is the angular velocity of the body
expressed in the body fixed frame, and $S(\cdot):\Re^3\mapsto \so$
is a skew mapping defined such that $S(x)y=x\times y$ for all
$x,y\in\Re^3$. $M\in\Re^3$ is the moment due to the potential. The
moment is determined by 
$S(M)=\deriv{U}{C}\Tp R-C\Tp \deriv{U}{C}$, or more explicitly,
\begin{gather}
M=r_1\times v_{r_1} + r_2\times v_{r_2} +r_3\times v_{r_3},
\end{gather}
where $r_i,v_{r_i}\in\Re^{1\times 3}$ are the $i$th row vectors of
$C$ and $\deriv{U}{C}$, respectively. The derivation of the above 
equations can be found in~\cite{cca05}.

General numerical integration methods, including the popular
Runge-Kutta schemes, typically preserve neither first integrals nor
the characteristics of the configuration space, $\SO$. In particular, 
the orthogonal structure of the rotation matrices is not preserved 
numerically. 
To resolve these problems, a Lie group variational integrator for
the attitude dynamics of a rigid body is proposed
in~\cite{cca05}. This Lie group variational integrator is
described by the discrete time equations.
\begin{gather}
h S(J\omega_k+\frac{h}{2} M_k) = F_k J_d - J_dF_k\Tp,\label{eqn:findf0}\\
C_{k+1} = C_k F_k,\label{eqn:updateR0}\\
J\omega_{k+1} = F_k\Tp J\omega_k +\frac{h}{2} F_k\Tp M_k
+\frac{h}{2}M_{k+1},\label{eqn:updatew0}
\end{gather}
where $J_d\in\Re^3$ is a nonstandard moment of inertia matrix
defined by $J_d=\frac{1}{2}\tr{J}I_{3\times 3}-J$, and $F_k\in\SO$
is the relative attitude over an integration step. The constant
$h\in\Re$ is the integration step size, and the subscript $k$
denotes the $k$th integration step. This integrator yields a map
$(C_k,\omega_k)\mapsto(C_{k+1},\omega_{k+1})$ by solving
\refeqn{findf0} to obtain $F_k\in\SO$ and substituting it into
\refeqn{updateR0} and \refeqn{updatew0} to obtain $C_{k+1}$ and
$\omega_{k+1}$.

Since this integrator does not use a local parameterization, the
attitude is defined globally without singularities. It preserves the
orthogonal structure of $\SO$ because the rotation matrix is updated
by a multiplication of two rotation matrices in \refeqn{updateR0}. 
This integrator is obtained from a discrete variational principle, 
and it exhibits the characteristic symplectic and momentum preservation 
properties, and good energy behavior characteristic of variational 
integrators. We use \refeqn{findf0}, \refeqn{updateR0}, and 
\refeqn{updatew0} in the following development of the attitude estimator.

\subsection{Uncertainty Ellipsoid}

The configuration space of the attitude dynamics is $\SO$, so the
state evolves in $\T\SO$. Thus the corresponding uncertainty ellipsoid 
is a submanifold of $\T\SO$. An uncertainty ellipsoid centered at 
$(\hat{C},\hat\omega)$ is induced from an uncertainty ellipsoid in $\Re^6$, 
using the Lie algebra $\so$;
\begin{align}
    &\mathcal{E}(\hat{C},\hat{\omega},P)
    \nn \\ &= \braces{C\in\SO,\,\omega\in\Re^3 \,\Big|\,
    \bracket{\zeta\Tp,\,\delta\omega\Tp} P^{-1}\begin{bmatrix}\zeta\\\delta\omega
    \end{bmatrix}\leq 1},\nonumber\\
        &= \braces{C\in\SO,\,\omega\in\Re^3 \,\Big|\,
    \begin{bmatrix}\zeta\\\delta\omega\end{bmatrix}\in\mathcal{E}_{\Re^6}
	(0_{6},P)},\label{eqn:ueso}
\end{align}
where $S(\zeta)=\mathrm{logm} \parenth{\hat{C}\Tp C}\in\so$,
$\delta\omega=\omega-\hat{\omega}\in\Re^3$, and $P\in\Re^{6\times
6}$ is a symmetric positive definite matrix. Equivalently, an
element $(C,\omega)\in\mathcal{E}(\hat{C},\hat{\omega},P)$ can be
written as
\begin{align*}
    C & = \hat{C} e^{ S(\zeta)},\\
    \omega & = \hat{\omega} + \delta \omega,
\end{align*}
for $x=\bracket{\zeta\Tp,\,\delta\omega\Tp}\Tp\in\Re^6$ satisfying
$x\Tp P^{-1}x\leq 1$.
%

\subsection{Measurement error model}
We give the measurement error models for the direction vector and
for the angular velocity. The direction vector $b^i\in\S^2$ is
in the body fixed frame, and let $\tilde b^i\in\S^2$ denote the 
corresponding measured directions. Since we only measure directions, 
it is inappropriate to express the measurement error by a vector 
difference. Instead, we model it by rotation of the measured direction;
\begin{align}
{b}^i& = e^{S(\nu^i)} \tilde b^i,\nonumber\\
& \simeq \tilde b^i + S(\nu^i)\tilde b^i,\label{eqn:bi}
\end{align}
where $\nu^i\in\Re^3$ is the sensor error, which represents the
Euler axis of rotation vector from $\tilde b^i$ to $b^i$, and
$\norm{\nu^i}$ is the corresponding rotation angle in radians. The 
second equality assumes small measurement errors. The angular velocity 
measurement errors are modeled as
\begin{align}
\omega_k=\tilde\omega_k + \upsilon_k,\label{eqn:Omega}
\end{align}
where $\tilde\omega_k \in\Re^3$ is the measured angular velocity,
and $\upsilon_k\in\Re^3$ is an additive error.

We assume that the initial conditions and the sensor noise are
bounded by prescribed uncertainty ellipsoids.
\begin{gather}
    (C_0,\omega_0)\in\mathcal{E}(\hat{C}_0,\hat{\omega}_0,P_0),\label{eqn:P0}\\
    \nu_k^i\in\mathcal{E}_{\Re^3}(0,S^i_k),\label{eqn:Sk}\\
    \upsilon_k\in\mathcal{E}_{\Re^3}(0,T_k)\label{eqn:Tk},
\end{gather}
where $P_0\in\Re^{6\times 6}$, $S_k^i, T_k \in\Re^{3\times 3}$ are
symmetric positive definite matrices that define the shape and the
size of the uncertainty ellipsoids. 

\section{Attitude Estimation with Angular Velocity Sensor}
In this section, we develop a deterministic estimator for the
attitude and the angular velocity of a rigid body assuming that both
attitude measurement and angular velocity measurements are
available. The estimator consists of three stages; flow propagation,
measurement, and filtered update. The propagation is to predict the
uncertainty ellipsoid in the future. The measurement is to find an 
uncertainty ellipsoid in the state space using the measurements and 
the measurement error model. The filtered update finds a new estimate 
using the predicted uncertainty ellipsoid and the measured uncertainty 
ellipsoid.

The subscript $k$ denotes the $k$-th discrete index. This may not 
coincide with measurement instants as we may resolve the evolution 
of the trajectory more frequently than the frequency of the 
measurements. This enables us to deal with measurements that are 
rather infrequent, with nontrivial attitude evolution between the 
measurements. The superscript $f$ denotes the variables related to 
the flow update, and the superscript $m$ denotes the variables 
related to the measurement update. $\tilde\cdot$ denotes a variable 
measured by sensors, and $\hat\cdot$ denotes an estimated variable.

\subsection{Flow propagation}
Suppose that the attitude and the angular momentum at the $k$th step
lie in a given uncertainty ellipsoid:
\begin{align*}
    (C_k,\omega_k)\in\mathcal{E}(\hat{C}_k,\hat{\omega}_k,P_k),
\end{align*}
and suppose that new measurements are taken at the $k+l$th time step. 

The flow update obtains the the uncertainty ellipsoid at the $k+l$th 
step using the given uncertainty ellipsoid at the $k$th step. 
We assume that the given uncertainty ellipsoid at the $k$th step is 
sufficiently small that the states in the uncertainty ellipsoid can 
be approximated by linearized equations of motion. This guarantees 
that the boundary of the state uncertainties at the $k+l$th step 
remains an ellipsoid. 

\textit{Center:} For the given center, $(\hat{C}_k,\hat{\omega}_k)$, 
the center of the uncertainty ellipsoid $(\hat{C}_{k+1}^{f},
\hat{\omega}_{k+1}^{f})$ is obtained from the discrete equations of motion, 
\refeqn{findf0}, \refeqn{updateR0}, and \refeqn{updatew0}:
\begin{gather}
h S(J\hat{\omega}_k+\frac{h}{2} \hat{M}_k) = \hat{F}_k J_d - J_d
\hat{F}_k\Tp,\label{eqn:findf}\\
\hat{C}_{k+1}^{f} = \hat{C}_k \hat{F}_k,\label{eqn:updateR}\\
J\hat{\omega}_{k+1}^{f} = \hat{F_k}\Tp
\hat{\omega}_k+\frac{h}{2}\hat{F_k}\Tp \hat{M_k} +\frac{h}{2}
\hat{M}_{k+1}.\label{eqn:updatePi}
\end{gather}
This integrator yields a map $(\hat C_k,\hat\omega_k)\mapsto(\hat
C^f_{k+1},\hat\omega_{k+1}^f)$, and this process can be repeated to
find the center at the $k+l$th step, $(\hat
C^f_{k+l},\hat\omega_{k+l}^f)$.

\textit{Uncertainty matrix:} The uncertainty matrix is obtained by 
linearizing the above discrete equations of motion. At the $(k+1)$th step, 
the state is given by perturbations from the center $(\hat C^f_{k+l},
\hat\omega_{k+l}^f)$ as
\begin{align*}
    C_{k+1}&=\hat{C}_{k+1}^{f} e^{S(\zeta_{k+1}^{f})},\\
    \omega_{k+1}&=\hat{\omega}_{k+1}^{f}+\delta\omega_{k+1}^{f},
\end{align*}
for some $\zeta_{k+1}^{f},\delta\omega_{k+1}^{f}\in\Re^3$. 
Assume that the uncertainty ellipsoid at the $k$th step is sufficiently
small. Then, $\zeta_{k+1}^{f},\delta\omega_{k+1}^{f}$ are given by the 
following linear equations in~\cite{acc06}:
\begin{align*}
\begin{bmatrix}\zeta_{k+1}^{f}\\\delta\omega_{k+1}^{f}\end{bmatrix}
& = \begin{bmatrix} \mathcal{A}_k^f & \mathcal{B}_k^f \\
\mathcal{C}_k^f & \mathcal{D}_k^f \end{bmatrix}
\begin{bmatrix}\zeta_{k}\\\delta\omega_{k}\end{bmatrix},
\end{align*}
where $\mathcal{A}^f_k, \mathcal{B}^f_k, \mathcal{C}^f_k,
\mathcal{D}^f_k,\in\Re^{3\times 3}$ can be suitably defined.
Equivalently, we rewrite the above equation as
\begin{align*}
x_{k+1}^{f} & = A_k^f x_k,
\end{align*}
where $x_k=[\zeta_k^T,\delta\omega_k^T]^T\in\Re^6$,
$A_k^f\in\Re^{6\times 6}$. Since 
$(C_k,\omega_k)\in\mathcal{E}(\hat{C}_k,\hat{\omega}_k,P_k)$,
$x_k\in\mathcal{E}_{\Re^6}(0,P_k)$ by the definition of the
uncertainty ellipsoid given in \refeqn{ueso}, we can show that
\begin{align*}
A_k^f x_k&\in\mathcal{E}_{\Re^6}\!\parenth{0,A_k^f P_k
\parenth{A_k^f}\Tp}.
\end{align*}
Thus, the uncertainty matrix at the $k+1$th step is given by
\begin{align}
P_{k+1}^f & = A_k^f P_k \parenth{A_k^f}\Tp.\label{eqn:Pkpf}
\end{align}
The above equation can be applied repeatedly to find the uncertainty
matrix at the $k+l$th step.
In summary, the uncertainty ellipsoid at the ($k+l$)th
step is computed using \refeqn{findf}, \refeqn{updateR},
\refeqn{updatePi}, and \refeqn{Pkpf} as:
\begin{align}\label{eqn:flow}
    (C_{k+l},\omega_{k+l})\in\mathcal{E}(\wh{C}_{k+l}^{f},
\wh{\omega}_{k+l}^{f},P_{k+l}^f),
\end{align}

\subsection{Measurement update}
The measured attitude and angular velocity have uncertainties due to 
sensor errors. However, we can find a uncertainty bound on the 
states because we assume that the sensor errors are bounded by known 
uncertainty ellipsoids. The measurement
update obtains an uncertainty ellipsoid in the state space using
the measurements and the sensor error models. 

\textit{Center:} The center of the uncertainty ellipsoid,
$(\wh{C}_{k+l}^{m},\wh\omega_{k+l}^m)$ is obtained from the 
measurements. The attitude is determined by measuring the directions
to the known points in the inertial frame. Let the measured
directions to the known points be $\wt B_{k+l}=\bracket{\wt
b^1,\wt b^2,\cdots,\wt b^m}\in\Re^{3\times m}$. Then, the
attitude $\wh{C}_{k+l}^{m}$ satisfies the following necessary and
sufficient condition given in \refeqn{san2006}.
\begin{gather}
\parenth{\wh{C}_{k+l}^m}\Tp\wt L_{k+l}-\wt L_{k+l}\Tp\wh{C}_{k+l}^m=0,
\label{eqn:meaR}
\end{gather}
where $\wt L_{k+l}=E_{k+l} W_{k+l}\wt{B}_{k+l}\Tp\in\Re^{3\times 3}$. 
The solution of \refeqn{meaR} is obtained by a QR factorization of 
$\wt L_{k+l}$ as given in Theorem \ref{atde}.
\begin{align}
\wh{C}_{k+l}^m=\parenth{Q\sqrt{(RR\Tp)^{-1}}Q\Tp} \wt L_{k+l},
\label{eqn:meaRex}
\end{align}
where $Q\in\SO$ and $R\in\Re^{3\times 3}$ is upper triangular such 
that $\wt L_{k+l}=QR$. The angular velocity is measured directly 
by sensors;
\begin{align}
\wh{\omega}_{k+l}^{m}= \wt\omega_{k+l}.\label{eqn:meaw}
\end{align}

\textit{Uncertainty matrix:} We can represent the actual state at
the $k+l$th step as follows:
\begin{align}
    C_{k+l}&=\wh{C}_{k+l}^{m} e^{S(\zeta_{k+l}^{m})},\label{eqn:Rkpm}\\
    \omega_{k+l}&=\wh{\omega}_{k+l}^{m}+\delta\omega_{k+l}^{m},\label{eqn:Omegakpm}
\end{align}
for $\zeta_{k+l}^{m},\delta\omega_{k+l}^{m}\in\Re^3$. The uncertainty 
matrix is obtained by finding an ellipsoidal bound for
$\zeta_{k+l}^{m},\delta\omega_{k+l}^{m}$. 

For the attitude, we transform the uncertainties in the directional sensors 
into the uncertainties in the rotation matrix by \refeqn{meaR}. 
The actual matrix of body direction vectors ${B}_{k+l}$ and the actual 
attitude $C_{k+l}$ also satisfy \refeqn{meaRex};
\begin{gather}
C_{k+l}\Tp L_{k+l}- L_{k+l}\Tp {C}_{k+l}=0,\label{eqn:meaRac}
\end{gather}
where $L_{k+l}=E_{k+l} W_{k+l}{B}_{k+l}\Tp\in\Re^{3\times 3}$. 
Using the identity, $S(x)A+A\Tp S(x)=S(\braces{\tr{A}I_{3\times 3}-A}x)$ 
for $A\in\Re^{3\times 3},x\in\Re^3$, \refeqn{meaRac} can be written 
in the vector form
\begin{align*}
\braces{\tr{\parenth{\wh{C}_{k+l}^m}\Tp \wt{L}_{k+l}}-\parenth{
\wh{C}_{k+l}^m}\Tp \wt{L}_{k+l}}\zeta_{k+l}^m \nn \\ =-\sum_{i=1}^m 
w_i \braces{\tr{\tilde{b}^i (e^i)\Tp\hat{C}_{k+l}^m}I_{3\times 3}-
\tilde{b}^i (e^i)\Tp\hat{C}_{k+l}^m}\nu^i.
\end{align*}
Then, we obtain
\begin{align}\label{eqn:zetakpm}
\zeta_{k+l}^m &= \sum_{i=1}^m \mathcal{A}_{k+l}^{m,i}\nu^i,
\end{align}
where
\begin{align}
\mathcal{A}_{k+l}^{m,i} & = -\braces{\tr{\parenth{\hat{C}_{k+l}^m}\Tp
\tilde{L}_{k+l}}-\parenth{\hat{C}_{k+l}^m}\Tp \tilde{L}_{k+l}}^{-1}
\nn \\ & w_i \braces{\tr{\tilde{b}^i (e^i)\Tp
\hat{C}_{k+l}^m}I_{3\times 3}-\tilde{b}^i (e^i)\Tp \hat{C}_{k+l}^m}.
\end{align}
This equation expresses the error in the measured attitude as a 
linear combination of the directional sensor errors. 

The perturbation of the angular velocity $\delta\omega_{k+l}^{m}$ 
is equal to the angular velocity measurement error $\upsilon_{k+l}$. 
Substituting \refeqn{Omegakpm} into \refeqn{Omega}, we obtain
\begin{align}\label{eqn:delPikpm}
\delta\omega_{k+l}^m = \upsilon_{k+l}.
\end{align}

Define the error states 
$x_{k+l}^m=\bracket{\parenth{\zeta^m_{k+l}}\Tp,\,\parenth{\delta
\omega^m_{k+l}}\Tp}\Tp\in\Re^6$.
Using \refeqn{zetakpm} and \refeqn{delPikpm},
\begin{align*}
x_{k+l}^m & = H_1 \sum_{i=1}^m \mathcal{A}_{k+l}^{m,i} \nu^i_{k+l} +
H_2 \upsilon_{k+l},
\end{align*}
where $H_1=[I_{3\times 3},\, 0_{3\times 3}]\Tp,H_2=[0_{3\times 3},\,
I_{3\times 3}]\Tp\in\Re^{6\times 3}$ which expresses $x_{k+l}^m$ as
a linear combination of the sensor errors $\nu^i$ and $\upsilon$.
From \refeqn{Sk} and \refeqn{Tk}, each term on the right hand side 
is in the following uncertainty ellipsoids:
\begin{align*}
H_1 \mathcal{A}_{k+l}^{m,i} \nu^i_{k+l} & \in
\mathcal{E}_{\Re^6}\parenth{0,H_1 \mathcal{A}_{k+l}^{m,i}S^i_{k+l}
\parenth{\mathcal{A}_{k+l}^{m,i}}\Tp H_1\Tp},\\
H_2 \upsilon_{k+l} & \in \mathcal{E}_{\Re^6}\parenth{0,H_2 T_{k+l}
 H_2\Tp}.
\end{align*}
The measurement update finds a minimal ellipsoid containing the 
vector sum of these uncertainty ellipsoids. Expressions for a minimal 
ellipsoid containing multiple ellipsoids are given in~\cite{MaNo96} 
and~\cite{DuWaPo01}, and $P_{k+l}^m$ is given by
{\small 
\begin{align}
& P_{k+l}^m  = \left\{\sum_{i=1}^m \sqrt{\tr{H_1
\mathcal{A}_{k+l}^{m,i}S^i_{k+l}\parenth{\mathcal{A}_{k+l}^{m,i}}\Tp H_1\Tp}}
\right. \nn \\ & \left. +\sqrt{\tr{H_2 T_{k+l} H_2\Tp}}\right\} \left\{
\sum_{i=1}^m \frac{H_1\mathcal{A}_{k+l}^{m,i}S^i_{k+l}
\parenth{\mathcal{A}_{k+l}^{m,i}}\Tp H_1\Tp}{\sqrt{\tr{H_1\mathcal{A}_{k+l}^{m,i}
S^i_{k+l}\parenth{\mathcal{A}_{k+l}^{m,i}}\Tp H_1\Tp}}}\right. \nn \\ & \left.
+\frac{H_2 T_{k+l} H_2\Tp}{\sqrt{\tr{H_2 T_{k+l}H_2\Tp}}}\right\}.\label{eqn:Pmkl}
\end{align}}

In summary, the measured uncertainty ellipsoid at the $k+l$th step
is defined by \refeqn{meaRex}, \refeqn{meaw}, and \refeqn{Pmkl};
\begin{align}\label{eqn:mea}
    (C_{k+l},\omega_{k+l})\in\mathcal{E}(\wh{C}_{k+l}^{m},\wh{\omega}_{k+l}^{m},
	P_{k+l}^m).
\end{align}

\subsection{Filtering procedure}
The filtering procedure is to find a new uncertainty ellipsoid
compatible with the predicted and the measured uncertainty ellipsoids. 
From \refeqn{flow} and \refeqn{mea}, the state at $k+l$th step lies 
in the intersection
\begin{align}\label{eqn:filter}
    (C_{k+l},\omega_{k+l})\in\mathcal{E}(\hat{C}_{k+l}^{f},
	\hat{\omega}_{k+l}^{f},P_{k+l}^f)\bigcap \nn \\
    \mathcal{E}(\hat{C}_{k+l}^{m},\hat{\omega}_{k+l}^{m},P_{k+l}^m).
\end{align}
Since it is inefficient to describe an irregular subset like the 
intersection of two ellipsoids in the state space numerically, we 
find a minimal uncertainty ellipsoid containing the intersection.
We omit the subscript $(k+l)$ in this subsection for convenience.

The measurement uncertainty ellipsoid, $\mathcal{E}(\hat{C}^{m},
\hat{\omega}^{m},P^m)$, is identified by its center $(\hat{C}^{m},
\hat{\omega}^{m})$, and the uncertainty ellipsoid in $\Re^6$:
\begin{align}
(\zeta^m,\delta\omega^m)\in\mathcal{E}_{\Re^6}(0_{6\times
1},P^m),\label{eqn:ellRm}
\end{align}
where $S(\zeta^m)=\mathrm{logm} \parenth{\hat{C}^{m,T} C}\in\so$,
$\delta\omega^m=\omega-\hat{\omega}^m\in\Re^3$. Similarly, the
predicted uncertainty ellipsoid, $\mathcal{E}(\hat{C}^{f},
\hat{\omega}^{f},P^f)$, is identified by its center $(\hat{C}^{f},
\hat{\omega}^{f})$, and the uncertainty ellipsoid in $\Re^6$:
\begin{align}
(\zeta^f,\delta\omega^f)\in\mathcal{E}_{\Re^6}(0_{6\times
1},P^f),\label{eqn:ellRf}
\end{align}
where $S(\zeta^f)=\mathrm{logm} \parenth{\hat{C}^{f,T} C}\in\so$,
$\delta\omega^f=\omega-\hat{\omega}^f\in\Re^3$. 

Define $\hat\zeta^{mf},\delta\hat\omega^{mf}\in\Re^3$ such that
\begin{align}
\hat{C}^f&=\hat{C}^m e^{S(\hat\zeta^{mf})},\label{eqn:errz}\\
\hat{\omega}^f& =\hat{\omega}^m
+\delta\hat\omega^{mf}.\label{eqn:errPi}
\end{align}
Thus, $\hat\zeta^{mf},\delta\hat\omega^{mf}$ gives the
difference between the centers of the two ellipsoids.
Using \refeqn{errz} and \refeqn{errPi} we get
\begin{align}
C^f & = \hat{C}^m e^{S(\hat\zeta^{mf})} e^{S(\zeta^f)},\nonumber\\
& \simeq \hat{C}^m e^{S(\hat\zeta^{mf}+\zeta^f)},\\
\omega^f & = \hat{\omega}^m + \parenth{\delta\hat\omega^{mf} +
\delta{\omega}^f},
\end{align}
where we assumed that $\hat\zeta^{mf}, \zeta^f$ are sufficiently
small. Thus, the uncertainty ellipsoid obtained by the flow update,
$\mathcal{E}(\hat{C}^{f},\hat{\omega}^{f},P^f)$ is given by the
center $(\hat{C}^m,\hat{\omega}^m)$ of the measurement uncertainty 
ellipsoid and 
\begin{align}
\mathcal{E}_{\Re^6}( \hat{x}^{mf} ,P^f),\label{eqn:ellRmf}
\end{align}
where
$\hat{x}^{mf}=\bracket{\parenth{\hat\zeta^{mf}}\Tp,
\parenth{\delta\hat\omega^{mf}}\Tp}\Tp\in\Re^6$.

We seek a minimal ellipsoid that contains the intersection of two
uncertainty ellipsoids in $\Re^6$:
\begin{align}
\mathcal{E}_{\Re^6}(0_{6\times 1},P^m)\bigcap
\mathcal{E}_{\Re^6}(\hat{x}^{mf} ,P^f)\subset\mathcal{E}_{\Re^6}(\hat{x},P),
\end{align}
where $\hat{x}=[\hat\zeta\Tp,\delta\hat\omega\Tp]\Tp\in\Re^6$. We obtain 
$\hat{x}$ and $P$ as
\begin{align*}
\hat{x}&=L\hat{x}^{mf},\\
P&=\beta(q) (I-L)P^m,
\end{align*}
where
\begin{align*}
\beta(q) & = 1 + q - (\hat{x}^{mf})\Tp (P^{m})^{-1} L \hat{x}^{mf},\\
L & = P^{m} (P^{m} + q^{-1} P^f)^{-1}.
\end{align*}
The constant $q$ is chosen such that $\tr{P}$ is minimized. We
convert $\hat{x}$ to points in $\T\SO$ using the common center
$(\hat{C}^{m},\hat{\omega}^{m})$.

In summary, a new uncertainty ellipsoid at the $k+l$th step is
defined by
\begin{align}
    (C_{k+l},\omega_{k+l})\in\mathcal{E}(\hat{C}_{k+l},\hat{\omega}_{k+l},P_{k+l}),
\end{align}
where
\begin{align}
\hat{C}_{k+l}&=\hat{C}_{k+l}^m e^{S(\hat\zeta)},\\
\hat{\omega}_{k+l}& = \hat{\omega}_{k+l}^m + \delta\hat\omega,\\
P_{k+l}& = P.
\end{align}

The entire estimation procedure is repeated. The new uncertainty
ellipsoid is used to predict the uncertainty ellipsoid till the next
measurements are available, and the measurement update and the
filtering procedures are performed. 
The center of the new uncertainty ellipsoid provides point
estimates of the attitude and the angular velocity at the $k+l$th
step. The uncertainty matrix represents the ellipsoidal bound on 
uncertainty. The size of the uncertainty matrix characterizes the
accuracy of the estimates. If the size is small, we conclude that 
the estimates are accurate. This estimation scheme is optimal since 
the size of the new uncertainty ellipsoid is minimized. 
The eigenvector of the uncertainty matrix corresponding to the maximum
eigenvalue shows the direction of the maximum uncertainty.

\section{Numerical Simulation}
Numerical simulation results are given for the estimation of the
attitude dynamics of an uncontrolled rigid spacecraft in a circular
orbit about a large central body, including gravity gradient
effects. The detailed description of the on orbit spacecraft model
is presented in~\cite{acc06}.

The inertia of the spacecraft is chosen as
$\overline{J}=\mathrm{diag}\bracket{1,\,2.8,\,2}$, where overlines
denote normalized variables. The maneuver is an arbitrary large
attitude change completed in a quarter of the orbit,
$\overline{T}_f=\frac{\pi}{2}\,\overline{s}$. The initial conditions 
are chosen as
\begin{alignat*}{2}
C_0 &= \mathrm{diag}[-1,\,-1,\,1],&\;
\overline\omega_0&=[2.3160,\,0.4468,\,-0.5910]\Tp,\\
\hat C_0 &= I_{3\times 3},&\;
\hat{\overline\omega}_0&=[2.1160,\,0.5468,\,-0.8910]\Tp.
\end{alignat*}
The corresponding initial estimation errors are
$\norm{\zeta_0}=180\,\mathrm{deg}$,
$\norm{\delta\omega_0}=21.43\frac{\pi}{180}\,\mathrm{rad/\overline{s}}$.
The initial uncertainty matrix is given by
\begin{align*}
P_0 =
2\,\mathrm{diag}\bracket{\parenth{180\frac{\pi}{180}}^2[1,\,1,\,1],\,
\parenth{30\frac{\pi}{180}}^2[1,\,1,\,1]},
\end{align*}
so that $x_0^TP_0^{-1}x_0=0.7553\leq 1$.

\begin{figure*}[th]
    \centerline{\subfigure[Estimation error $\norm{\zeta_k}$, $\norm{\delta\Omega_k}$]{
    \includegraphics[width=0.45\textwidth]{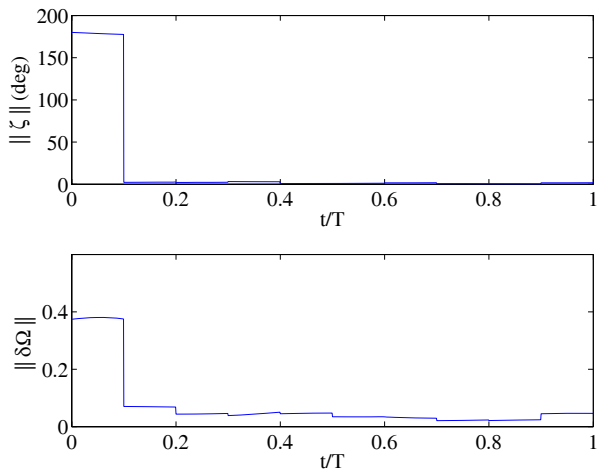}}
    \hfill
    \subfigure[Magnitude of Uncertainty $\tr{P_k}$]{
    \includegraphics[width=0.45\textwidth]{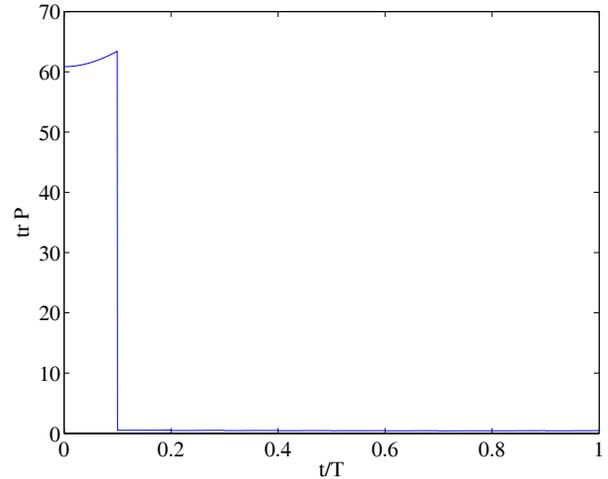}}
    }
    \caption{Attitude and angular velocity estimation errors with 
    measurements}\label{fig:full}
\end{figure*}

We assume that measurements are available ten times in a quarter orbit. 
The measurement uncertainty matrices are given by
\begin{align*}
S^i_k=\parenth{7\frac{\pi}{180}}^2 I_{3\times
3}\,\mathrm{rad^2},\; T_k=\parenth{7\frac{\pi}{180}}^2 I_{3\times
3}\,\mathrm{rad^2/\overline{s}^2}.
\end{align*}
\reffig{full} shows simulation results for a typical realization of the bounded uncertainties, where the plot on the 
left shows the attitude estimation error and the angular
velocity estimation error, and the right plot shows the size of
the uncertainty ellipsoid. The estimation errors and the size of
uncertainty decrease fast after the first estimation. The
terminal attitude error is less than $1\,\mathrm{deg}$.

\section{Conclusions}

The attitude estimation scheme presented here has no singularities since 
the attitude is represented by a rotation matrix, and the structure of 
the group of rotation matrices is preserved since it is updated by group 
operations in $\SO$ using the Lie group variational integrator. The 
attitude estimator is also robust to the distribution of the uncertainty 
and the sensor noise, since it is based on deterministic ellipsoidal 
bounds on the uncertainty. The effects of process noise can be included 
by modifying the prediction procedure.

Although not presented in this paper, we have obtained results for the 
modification of this scheme to the case when angular velocity measurements 
are not available. We intend to extend this estimation scheme to the 
combined attitude control and estimation problem for a rigid body in 
an attitude dependent potential, with the inclusion of process noise 
or disturbance forces. These topics will be dealt with in a future 
journal paper.

\section{ACKNOWLEDGMENTS}
The research of ML has been supported in part by NSF under Grant DMS-0504747.
The research of NHM has been supported in part by NSF under Grant ECS-0140053.

\end{document}